\documentclass[twoside,11pt]{article}

\usepackage{jmlr2e}
\usepackage{lineno}
\usepackage{hyperref}       
\usepackage{amsfonts}       
\usepackage{bbm}
\usepackage{amssymb}
\usepackage{amsmath}

\newtheorem{prop}{Proposition}

\firstpageno{1}

\begin{document}

\title{On Concentration Inequalities for Vector-Valued Lipschitz Functions}

\author{\name Dimitrios Katselis \email katselis@illinois.edu \\
       \addr ECE Department, University of Illinois at Urbana-Champaign, USA
       \AND
       \name Xiaotian Xie \email xx5@illinois.edu \\
       \addr ISE Department and Coordinated Science Lab, University of Illinois at Urbana-Champaign, USA
       \AND
       \name Carolyn L. Beck \email beck3@illinois.edu \\
       \addr ISE Department and Coordinated Science Lab, University of Illinois at Urbana-Champaign, USA
       \AND
       \name R. Srikant \email rsrikant@illinois.edu \\
       \addr ECE Department and Coordinated Science Lab, University of Illinois at Urbana-Champaign, USA
       }

\maketitle

\begin{abstract}
We derive two upper bounds for the probability of deviation of a vector-valued Lipschitz function of a collection of random variables from its expected value.  The resulting upper bounds can be tighter than bounds obtained by a direct application of a classical theorem due to Bobkov and G\"{o}tze. 
\end{abstract}


\section{Introduction}

In many statistical settings, vector-valued estimators naturally arise and determining their statistical rates is essential. As an example, assume that $X=(X_1,\ldots, X_n)$ corresponds to a random sample drawn from a product measure $\nu_{\theta}^{\otimes n}$, denoted by $X\sim \nu_{\theta}^{\otimes n}$, where $\theta$ is a parameter vector taking values in a compact set $\Theta \subset \mathbb{R}^k$ with a fixed dimension $k$ independent of $n$. In this case, $\hat{\theta}_n=f(X)$ is a candidate vector-valued estimator of $\theta$, where $f$ is a measurable function with respect to $X$.  In such settings, we are interested in knowing how close $f(X)$ is to its expectation $E_{\nu}[f(X)]$. Concentration inequalities for $f(X)$ when the deviation from  $E_{\nu}[f(X)]$ is measured in terms of a metric, often norm-induced, are important. 

The derivation of concentration bounds relies on imposing  smoothness conditions on $f$, which guarantee that $f$ is not very sensitive to any particular coordinate variable \cite{blm13, rs13, vh16}. This sensitivity is quantified either locally via gradients or globally via Lipschitz properties of $f$. Marton introduced the  \emph{transportation method} to establish concentration of measure for product measures and Markov chains \cite{m86, m96} by showing that transportation cost inequalities can be used to deduce concentration. In \cite{bg99}, Bobkov and G\"{o}tze extended Marton's argument into an \emph{equivalence} by showing that Wasserstein distances and relative entropies are comparable only when the moment-generating functions of \emph{real-valued} Lipschitz functions defined with respect to the underlying metric can be controlled and vice versa. The connection between moment-generating functions and the two aforementioned indices of closeness of probability measures is established via the Gibbs variational principle \cite{vh16},\cite{zd98} or alternatively, via the Donsker-Varadhan lemma \cite{r15}. 

A key example in showing the connection between concentration of measure and Lipschitz functions is McDiarmid's or \emph{bounded-difference} inequality \cite{mc97}, traditionally viewed as a result of the  \emph{martingale approach} in establishing concentration \cite{blm13, rs13, vh16}. Let each $X_i$ take values in a measurable space $\mathcal{X}_i$ and equip the product space $\mathbb{X}=\mathcal{X}_1\times \cdot \times \mathcal{X}_n$ with the weighted Hamming metric $d_c(x,y)=\sum_{i=1}^n c_i \mathbbm{1}\{x_i\neq y_i\}$, where $x=(x_1,\ldots, x_n), y=(y_1,\ldots, y_n) \in \mathbb{X}$. For $X \in \mathbb{X}$ with independent entries, $f(X)$ is a sub-Gaussian random variable with parameter $\sigma^2=\sum_{i=1}^n c_i^2/4$ for every $f:\mathbb{X}\rightarrow \mathbb{R}$ which is $1$-Lipschitz with respect to $d_c$ by McDiarmid's inequality. Motivated by this example, a question of interest is for which measures $\nu$ on a metric space $(\mathbb{X},d)$ such that $X\sim \nu$, the random variable $f(X)$ is $\sigma^2$-sub-Gaussian for every \emph{real-valued} $1$-Lipschitz function $f:\mathbb{X}\rightarrow \mathbb{R}$. The answer to this question is given by the aforementioned theorem of Bobkov and G\"{o}tze \cite{bg99}. In this paper, we focus on Lipschitz mappings $f:\mathbb{X}\rightarrow \mathbb{Y}$ between metric spaces $(\mathbb{X},d_{\mathbb{X}})$ and $(\mathbb{Y},d_{\mathbb{Y}})$, where $\mathbb{Y}\subseteq \mathbb{R}^k$ and $d_{\mathbb{Y}}$ is any $\ell_p$-metric for $p\geq 1$. In the spirit of \cite{m86,m96,bg99}, we prove a concentration inequality when a transportation cost inequality can be shown to hold. We then provide a simple stationary measure estimation example for  Markov chains demonstrating that the derived inequality gives better results than directly applying the theorem of Bobkov and G\"{o}tze or by combining the aforementioned theorem with Boole's inequality. An interesting observation regarding the role of the particular $\ell_p$-norm is also highlighted via this example.

\section{Preliminaries and Main Result}

Let $(\mathbb{X},d_{\mathbb{X}})$ be a Polish space and $\rho$ be a (Borel) probability measure on $(\mathbb{X},d_{\mathbb{X}})$. The triplet  $(\mathbb{X},d_{\mathbb{X}},\rho)$ defines a metric probability space in the sense of Gromov \cite{rs13}, \cite{g07}. Given two metric spaces $(\mathbb{X}, d_{\mathbb{X}})$ and $(\mathbb{Y}, d_{\mathbb{Y}})$, let $f:\mathbb{X}\rightarrow \mathbb{Y}$ be a Lipschitz mapping with Lipschitz constant $\|f\|_{\rm Lip}$, i.e., $d_{\mathbb{Y}}(f(x),f(\tilde{x}))\leq \|f\|_{\rm Lip} d_{\mathbb{X}}(x,\tilde{x}),$ $\forall x,\tilde{x}\in \mathbb{X}$. In the following, the set of all such mappings will be denoted by ${\rm Lip}(\mathbb{X},\mathbb{Y},d_{\mathbb{X}},d_{\mathbb{Y}})$. Moreover,  let $\mathbb{P}_{1}(\mathbb{X})$ denote the set of all probability measures $\rho$ on $\mathbb{X}$ such that $E_{\rho}[d_{\mathbb{X}}(X,x_0)]<\infty$ holds for an arbitrary (and therefore for all) $x_0\in \mathbb{X}$. The $L^1$ Wasserstein distance between $\rho,\tilde{\rho}\in \mathbb{P}_1(\mathbb{X})$ is defined as
\begin{equation}\label{eq:defin1}
W_{1}(\rho,\tilde{\rho})=\inf_{X\sim \rho,\tilde{X}\sim \tilde{\rho}}E[d_{\mathbb{X}}(X,\tilde{X})]=\inf_{\gamma\in \Pi(\rho,\tilde{\rho})}\int_{\mathbb{X}\times \mathbb{X}}d_{\mathbb{X}}(x,\tilde{x})\gamma(dx,d\tilde{x}).
\end{equation}
The infimum in the first part is taken over all jointly distributed pairs $(X,\tilde{X})$ on the product space $\mathbb{X}^2=\mathbb{X}\times \mathbb{X}$ with marginals $\rho,\tilde{\rho}$, respectively. In the last part, $\Pi(\rho,\tilde{\rho})$ denotes the set of all possible couplings of $\rho,\tilde{\rho}$. An optimal coupling $\gamma^{*}\in \Pi(\rho,\tilde{\rho})$ achieving the infimum exists \cite{rs13}, \cite{v08}. Additionally, a different measure gauging the dissimilarity between two probability measures $\rho,\tilde{\rho}$ is the relative entropy or Kullback-Leibler divergence 
$
D(\tilde{\rho}\| \rho)=E_{\tilde{\rho}}\left[\log \frac{d\tilde{\rho}}{d\rho}\right]=E_{\rho}\left[\frac{d\tilde{\rho}}{d\rho}\log \frac{d\tilde{\rho}}{d\rho}\right]
$
for $\tilde{\rho}\ll \rho$ and $D(\tilde{\rho}\| \rho)=\infty$ otherwise, where $d\tilde{\rho}/d\rho$ is the Radon-Nikodym derivative of $\tilde{\rho}$ with respect to $\rho$ and $\tilde{\rho}\ll \rho$ denotes that $\tilde{\rho}$ is absolutely continuous with respect to $\rho$.

The following theorem provides a concentration inequality for $f(X)$ when the deviation from $E_{\nu}[f(X)]$ is measured in terms of the $\ell_2$-metric.

\begin{theorem}\label{thm:MainThm}
Let $(\mathbb{X},d_{\mathbb{X}})$ and $(\mathbb{Y},d_{\mathbb{Y}})$ be two  Polish spaces, where $\mathbb{Y}\subseteq \mathbb{R}^k$ and $d_{\mathbb{Y}}(y,\tilde{y})=\|y-\tilde{y}\|_2$ with $\|\cdot\|_2$ being the Euclidean norm.  Let $X$ be a random variable taking values in $\mathbb{X}$ and assume that $X\sim \nu$, where $\nu$ is a probability measure on $(\mathbb{X},d_{\mathbb{X}})$. Then, the inequality $W_{1}(\mu,\nu)\leq \sqrt{2\sigma^2 D(\mu\|\nu)}, \forall \mu$ implies that $\forall \epsilon>0$, $\forall \varepsilon\in (0,1]$ and for any $f\in {\rm Lip}(\mathbb{X},\mathbb{Y},d_{\mathbb{X}},\|\cdot\|_2)$ such that $\|E_{\nu}[f]\|_\infty<\infty$,
\begin{align}\label{eq:MainThm_1}
P(\|f(X)-E_{\nu}[f(X)]\|_2\geq \epsilon)\leq \min\left\{\left(1+\frac{2}{\varepsilon}\right)^k e^{-\frac{\epsilon^2 (1-\varepsilon)^2}{2 \sigma^2\|f\|_{\rm Lip}^2}}, 2^{\frac{k}{2}}e^{-\frac{\epsilon^2}{4\sigma^2 \|f\|_{\rm Lip}^2}}\right\}.
\end{align}
\end{theorem}  

The proof of this theorem is provided in Section \ref{sec:ProofMain}.  This result can be straightforwardly extended to  any $\ell_p$-metric for $p\geq 1, p\neq 2$. Define 
\begin{align}\label{eq:MainThm_3}
\tau_p&=\sup_{y\in f(\mathbb{X})}\frac{\|y-E_{\nu}[f(X)]\|_p}{\|y-E_{\nu}[f(X)]\|_2}.
\end{align}
Here, $y$ is assumed different from $E_{\nu}[f(X)]$ if $E_{\nu}[f(X)]\in f(\mathbb{X})$, since for $y=E_{\nu}[f(X)]$ (if $E_{\nu}[f(X)]\in f(\mathbb{X})$) the inequality $\|y-E_{\nu}[f(X)]\|_p\leq \tau_p \|y-E_{\nu}[f(X)]\|_2$ trivially holds  for any $\tau_p>0$. Then, (\ref{eq:MainThm_1}) can be replaced by 
\begin{align}\label{eq:MainThm_4}
P(\|f(X)-E_{\nu}[f(X)]\|_p\geq \epsilon)\leq \min\left\{\left(1+\frac{2}{\varepsilon}\right)^k e^{-\frac{\epsilon^2 (1-\varepsilon)^2}{2\tau_p^2 \sigma^2\|f\|_{\rm Lip}^2}}, 2^{\frac{k}{2}}e^{-\frac{\epsilon^2}{4\sigma^2\tau_p^2 \|f\|_{\rm Lip}^2}}\right\}
\end{align}
due to $\{\|f(X)-E_{\nu}[f(X)]\|_p\geq \epsilon\}\subseteq \{\|f(X)-E_{\nu}[f(X)]\|_2\geq \epsilon/\tau_p\}$ by (\ref{eq:MainThm_3}).

\section{Proof of Theorem \ref{thm:MainThm} and Additional Results}
\label{sec:ProofMain}

\noindent\textbf{Proof of the first bound in (\ref{eq:MainThm_1})}: The proof of the first bound relies on covering arguments; see \cite{v18,w19,lc20} and references therein for results based on such arguments.  For some $\varepsilon$ in the interval $(0,1]$, let $\mathcal{N}(\varepsilon)$ be an $\varepsilon$-net of the unit Euclidean  sphere $\mathbb{S}^{k-1}$ in $\mathbb{R}^k$ with cardinality $\left|\mathcal{N}(\varepsilon)\right|\leq \left(1+\frac{2}{\varepsilon}\right)^k$ \cite{v18}. Additionally, by Exercise 4.4.2 in \cite{v18}, 

\begin{equation}\label{eq:MainThm_proof_1}
\left\|f(X)-E_{\nu}[f(X)]\right\|_2\leq \frac{1}{1-\varepsilon}\sup_{w \in \mathcal{N}(\varepsilon)}\langle w, f(X)-E_{\nu}[f(X)] \rangle\ \ \ \ \text{a.s.}
\end{equation}
By (\ref{eq:MainThm_proof_1}) and  by Boole's inequality (union bound) we have that
\begin{align}\label{eq:MainThm_proof_2}
P(\|f(X)-E_{\nu}[f(X)]\|_2\geq \epsilon)&\leq P\left(\sup_{w \in \mathcal{N}(\varepsilon)}\langle w, f(X)-E_{\nu}[f(X)] \rangle \geq \epsilon (1-\varepsilon)\right) \nonumber\\&\leq \sum_{w\in \mathcal{N}(\varepsilon)}P\left(\langle w, f(X)-E_{\nu}[f(X)] \rangle \geq \epsilon (1-\varepsilon)\right)\nonumber\\&\leq |\mathcal{N}(\varepsilon)| \inf_{\lambda>0} e^{-\lambda \epsilon (1-\varepsilon)}E_{\nu}\left[e^{\lambda  \langle w_{*}, f(X)-E_{\nu}[f(X)] \rangle }\right],
\end{align}
where the last inequality follows from the Chernoff bound \cite{rs13} and
\[
w_{*}=\arg\max_{w\in \mathcal{N}(\varepsilon)} P\left(\langle w, f(X)-E_{\nu}[f(X)] \rangle \geq \epsilon (1-\varepsilon)\right).
\] 
We now note that  the function $ \langle w_{*}, f(x)-E_{\nu}[f(X)] \rangle$ is $ \|f\|_{\rm Lip}$-Lipschitz by the  Cauchy–-Schwarz inequality and $ \langle w_{*}, f(X)-E_{\nu}[f(X)] \rangle$ is mean zero.
Assuming that $W_{1}(\mu,\nu)\leq \sqrt{2\sigma^2 D(\mu\|\nu)}, \forall \mu$ holds, an application of the theorem of Bobkov and G\"{o}tze  implies that 
\begin{align}\label{eq:MainThm_proof_3}
P(\|f(X)-E_{\nu}[f(X)]\|_2\geq \epsilon)&\leq |\mathcal{N}(\varepsilon)| \inf_{\lambda>0} e^{-\lambda \epsilon (1-\varepsilon)+\frac{\lambda^2\sigma^2 \|f\|_{\rm Lip}^2}{2}}.
\end{align}
The exponent in the right-hand side of (\ref{eq:MainThm_proof_3}) is minimized for 
$
\lambda_{*}=\frac{\epsilon (1-\varepsilon)}{\sigma^2\|f\|_{\rm Lip}^2}>0
$
leading to 
\begin{align}\label{eq:MainThm_proof_4}
P(\|f(X)-E_{\nu}[f(X)]\|_2\geq \epsilon)&\leq |\mathcal{N}(\varepsilon)| e^{-\frac{\epsilon^2 (1-\varepsilon)^2}{2\sigma^2\|f\|_{\rm Lip}^2}}.
\end{align}
By employing the bound on $|\mathcal{N}(\varepsilon)|$ the desired result follows.

\noindent \textbf{Proof of the second bound in (\ref{eq:MainThm_1})}: Using a similar argument as in  \cite{hkz12}, let $Z\sim \mathcal{N}(0,I_k)$ be a standard Gaussian random vector, which is independent of $X$. Recall that $E[e^{\langle Z,q \rangle}]=e^{\frac{\|q\|_2^2}{2}}, \forall q\in \mathbb{R}^k$. For any $\lambda \in \mathbb{R}$ we note that
\begin{align*}
&E\left[e^{\lambda \langle Z,f(X)-E_{\nu}[f]\rangle}\right]\geq\\&  E\left.\left[e^{\lambda \langle Z, f(X)-E_{\nu}[f] \rangle}\right|\|f(X)-E_{\nu}[f(X)]\|_2\geq \epsilon\right]P(\|f(X)-E_{\nu}[f(X)]\|_2\geq \epsilon)=\\& E_{\nu}\left.\left[E_Z\left[e^{\lambda \langle Z, f(X)-E_{\nu}[f] \rangle}\right]\right|\|f(X)-E_{\nu}[f(X)]\|_2\geq \epsilon\right]P(\|f(X)-E_{\nu}[f(X)]\|_2\geq \epsilon)=\\& E_{\nu}\left.\left[e^{\frac{\lambda^2 \|f(X)-E_{\nu}[f] \|_2^2}{2}}\right|\|f(X)-E_{\nu}[f(X)]\|_2\geq \epsilon\right]P(\|f(X)-E_{\nu}[f(X)]\|_2\geq \epsilon)\geq \\& e^{\frac{\lambda^2 \epsilon^2}{2}} P(\|f(X)-E_{\nu}[f(X)]\|_2\geq \epsilon)
\end{align*}
or 
\begin{align}\label{eq:MainThm_proof_5}
P(\|f(X)-E_{\nu}[f(X)]\|_2\geq \epsilon)\leq e^{-\frac{\lambda^2 \epsilon^2}{2}} E\left[e^{\lambda \langle Z,f(X)-E_{\nu}[f]\rangle}\right].
\end{align}
We now note that the function $ \langle Z, f(x)-E_{\nu}[f(X)] \rangle$ is $\|Z\|_2\|f\|_{\rm Lip}$-Lipschitz when $Z$ is fixed (by the  Cauchy–-Schwarz inequality) and $ \langle Z, f(X)-E_{\nu}[f(X)] \rangle$ is mean zero. 
Assuming that $W_{1}(\mu,\nu)\leq \sqrt{2\sigma^2 D(\mu\|\nu)}, \forall \mu$ holds, an application of the theorem of Bobkov and G\"{o}tze  implies that 
\begin{align*}
E\left[e^{\lambda \langle Z,f(X)-E_{\nu}[f]\rangle}\right]&=E_{Z}\left[E_{\nu}\left[e^{\lambda \langle Z,f(X)-E_{\nu}[f]\rangle}\right]\right] \leq E_Z\left[e^{\frac{\lambda^2 \|Z\|_2^2\sigma^2 \|f\|_{\rm Lip}^2}{2}}\right].  
\end{align*}
Moreover, $\|Z\|_2^2$ is a chi-squared random variable with $k$ degrees of freedom and for such a random variable 
$
E\left[e^{t\|Z\|_2^2}\right]=\frac{1}{(1-2t)^{k/2}}, \ \ t<\frac{1}{2}
$. Therefore, we conclude that
\begin{align}\label{eq:MainThm_proof_6}
E\left[e^{\lambda \langle Z,f(X)-E_{\nu}[f]\rangle}\right]\leq E_Z\left[e^{\frac{\lambda^2 \|Z\|_2^2\sigma^2 \|f\|_{\rm Lip}^2}{2}}\right]= \frac{1}{\left(1-\lambda^2 \sigma^2 \|f\|_{\rm Lip}^2\right)^{k/2}},  
\end{align}
where the last equality holds for any  $\lambda\in \mathbb{R}$ such that $\lambda^2 \sigma^2 \|f\|_{\rm Lip}^2<1$. By combining (\ref{eq:MainThm_proof_5}) and (\ref{eq:MainThm_proof_6}) we obtain 
\begin{align}\label{eq:MainThm_proof_7}
P(\|f(X)-E_{\nu}[f(X)]\|_2\geq \epsilon)\leq \frac{e^{-\frac{\lambda^2 \epsilon^2}{2}}}{\left(1-\lambda^2 \sigma^2 \|f\|_{\rm Lip}^2\right)^{k/2}}, \ \ |\lambda|<\frac{1}{\sigma \|f\|_{\rm Lip}}. 
\end{align}
Choosing $\lambda=\frac{1}{\sqrt{2}\sigma \|f\|_{\rm Lip}}$, we conclude that 
\begin{align}
P(\|f(X)-E_{\nu}[f(X)]\|_2\geq \epsilon)\leq 2^{\frac{k}{2}}e^{-\frac{\epsilon^2}{4\sigma^2 \|f\|_{\rm Lip}^2}}.
\end{align}
 \hfill $\square$
 
 In the previous proof and more specifically in (\ref{eq:MainThm_proof_3}) and (\ref{eq:MainThm_proof_6}), the theorem of Bobkov and G\"{o}tze for the mean zero, vector-valued function $f(X)-E_{\nu}[f(X)]$ has been applied via a real-valued function of the form $\langle g,f(X)-E_{\nu}[f(X)]\rangle$ with $g=w_{*}$ and $g=Z$, respectively. This suggests the following extension of the theorem of Bobkov and G\"{o}tze for vector-valued Lipschitz functions: 
 
 \begin{prop}\label{prop:Bobkov-Gotze_extension}
(\textbf{Theorem of Bobkov and G\"{o}tze for Vector-Valued Functions})
Let $(\mathbb{X},d_{\mathbb{X}})$ and $(\mathbb{Y},d_{\mathbb{Y}})$ be two  Polish spaces, where $\mathbb{Y}\subseteq \mathbb{R}^k$ and $d_{\mathbb{Y}}$ is an $\ell_p$-metric for some $p\geq 1$. Let $X$  be a random variable taking values in $\mathbb{X}$ and assume that $X\sim \nu$, where $\nu$ is a  probability measure on $(\mathbb{X},d_{\mathbb{X}})$. Then, the following statements are equivalent:
\begin{enumerate}
\item $W_{1}(\mu,\nu)\leq \sqrt{2\sigma^2 D(\mu\|\nu)}, \forall \mu$ 
\item $f(X)-E_{\nu}[f(X)]$, for every function $f\in {\rm Lip}(\mathbb{X},\mathbb{Y},d_{\mathbb{X}},d_{\mathbb{Y}})$ such that $\|E_{\nu}[f]\|_\infty<\infty$, is a sub-Gaussian vector with $\sigma^2\|h\|_q^2\|f\|_{\rm Lip}^2$-sub-Gaussian one-dimensional marginals $\langle h,f(X)-E_{\nu}[f(X)]\rangle$, i.e., 
 \begin{equation}\label{eq:Bobkov_Gotze_extension_1}
E_{\nu}\left[e^{\lambda \langle h,f(X)-E_{\nu}[f(X)]\rangle}\right]\leq e^{\frac{\lambda^2 \|h\|_q^2 \sigma^2 \|f\|_{\rm Lip}^2}{2}}, \forall h\in \mathbb{R}^{k}, \forall \lambda\in \mathbb{R}.
\end{equation}
 Here, $\|\cdot\|_q$ corresponds to the dual norm of $\|\cdot\|_p$. 
 \end{enumerate}
\end{prop}  
The direction ``\textbf{$1.$ implies $2.$}'' can be obtained by the theorem of Bobkov and G\"{o}tze for the real-valued function $\langle h, f(x)-E_{\nu}[f(X)]\rangle$ by invoking H\"{o}lder's inequality to show that the corresponding Lipschitz constant is at most $\|h\|_q \|f\|_{\rm Lip}$. The direction ``\textbf{$2.$ implies $1.$}'' is a direct consequence of the theorem of Bobkov and G\"{o}tze by choosing $f$ to have only one nonzero coordinate, e.g., $f=\tilde{f}e_1$ and $h=e_1$. Here, $\tilde{f}:\mathbb{X}\rightarrow \mathbb{R}$ is any real-valued $\|f\|_{\rm Lip}$-Lipschitz function and $e_1$ is the first element of the canonical basis in $\mathbb{R}^k$.

\section{Example}

 Consider a sample $X_{0:n-1}=\{X_0,X_1,\ldots, X_{n-1}\}$ of size $n$ drawn from an ergodic, discrete-time, finite-state Markov chain $(X_k)_{k\geq 0}$ with state space $\mathbb{E}=[K]=\{1,2,\ldots, K\}$, transition matrix $\mathrm{P}=[\mathrm{P}_{ij}]$ and $X_0\sim \varrho$, where $\varrho$ denotes the initial measure of the chain. We denote such a Markov chain by $(\mathrm{P},\varrho)$ and the corresponding stationary chain by $(\mathrm{P},\pi)$, where $\pi$ is the underlying invariant measure. Let the chain be $r$-contractive with Dobrushin coefficient  $r<1$. Consider the natural plug-in estimators for the stationary probabilities: 
\[
\hat{\pi}_i(X_{0:n-1})=\frac{1}{n}\sum_{k=0}^{n-1}\mathbbm{1}(X_k=i).
\] 
By the Ergodic Theorem for Markov chains \cite{b13}, $\hat{\pi}_i\rightarrow \pi_i, \forall i\in \mathbb{E}$ with probability $1$ as $n\rightarrow \infty$ for any initial measure $\varrho$. 

 The distance of $(\mathrm{P},\varrho)$ from stationarity can be quantified by the (nonstationarity) index \cite{p15}
\begin{equation}\label{eq:index}
\left\|\frac{\varrho}{\pi}\right\|^2_{2,\pi}=E_{\pi}\left[\left(\frac{d\varrho}{d\pi}\right)^2\right]=\sum_{i\in \mathbb{E}} \left[\frac{\varrho(i)}{\sqrt{\pi(i)}}\right]^2,
\end{equation} 
where the first equality corresponds to the general definition of the index for $\varrho\ll \pi$ and the second equality is the specialization of this definition to our setting. Furthermore, $1\leq \left\|\varrho/\pi\right\|_{2,\pi}\leq \infty$ and $\|\cdot\|_{2,\pi}$ is the norm induced by the inner product $\langle f,g \rangle_{\pi}=\sum_{i\in \mathbb{E}} f(i)g(i)\pi(i)$ in $\ell^2(\pi)$ \cite{lpw09}. Due to ergodicity, $\min_{i\in \mathbb{E}}\pi(i)>0$ and also $\left\|\varrho/\pi\right\|_{2,\pi}\leq 1/\sqrt{\min_{i\in \mathbb{E}}\pi(i)}$. Additionally, $\left\|\varrho/\pi\right\|_{2,\pi}=1$ for $\varrho=\pi$  and $\left\|\varrho/\pi\right\|_{2,\pi}=\infty$ if $\varrho$ is not absolutely continuous with respect to $\pi$.

The index in (\ref{eq:index}) is useful in our context due to the following theorem \cite{p15}:
\begin{theorem}\label{thm:Paulin}
 Let $X_{0:n-1}$ be a sample drawn from a time-homogeneous Markov chain $(\mathrm{P},\varrho)$ with state space $\mathbb{E}$ and stationary measure $\pi$. Then for any measurable function $g:\mathbb{E}^{n}\rightarrow \mathbb{R}$ and $\forall \epsilon>0$,
\[
P_{\varrho}\left(g(X_{0:n-1}\right)\geq \epsilon)\leq \left\|\frac{\varrho}{\pi}\right\|_{2,\pi}\sqrt{P_{\pi}(g(X_{0:n-1})\geq \epsilon)},
\]
where $P_{\varrho}$ is the law of $(\mathrm{P},\varrho)$ and  $P_{\pi}$ is the law of $(\mathrm{P},\pi)$. 
\end{theorem} 
Our goal is  to bound $P_{\varrho}(\|\hat{\pi}-\pi\|_{p}\geq \epsilon)=P_{\varrho}(\|\hat{\pi}-E_{\pi}[\hat{\pi}]\|_p\geq \epsilon)$ for  $p\geq 1$. To tackle the problem within the transportation method framework, we will use the following theorem due to Marton \cite{m96} adapted to our setting:

\begin{theorem}\label{thm:Marton} 
 Consider a Markov chain $(X_k)_{k\geq 0}$ with a finite state space $\mathcal{X}$, transition matrix $P=[P_{ij}]$ and Dobrushin coefficient $r<1$. 
For $x_{0:n-1},\tilde{x}_{0:n-1}\in \mathcal{X}^{n}$ let $
\tilde{d}_{1,\mathcal{X}^{n}}(x_{0:n-1},\tilde{x}_{0:n-1})=\sum_{k=0}^{n-1}\mathbbm{1}(x_k\neq \tilde{x}_k)$. 
Then, 
\[
W_{1}\left(\mu,\nu\right)\leq \left[\frac{n}{2(1-r)^2}D(\mu\| \nu)\right]^{1/2},
\]
where $\nu=P_{\varrho}$ is the measure on $(\mathbb{X}=\mathcal{X}^{n},d_{\mathbb{X}}=\tilde{d}_{1,\mathcal{X}^{n}})$ due to the Markov chain  starting at some arbitrary initial measure $\varrho$ and $\mu$ is any  measure on $(\mathbb{X},d_{\mathbb{X}})$.
\end{theorem}

We now compare different approaches for bounding $P_{\varrho}(\|\hat{\pi}-\pi\|_{p}\geq \epsilon)$ and show that the bound obtained in Theorem 1 gives better results than other, more direct applications of the theorem of Bobkov and G\"{o}tze.

\noindent\textbf{Approach 1: Direct application of the theorem of Bobkov and G\"{o}tze}. Let $f(x)=f(x_{0:n-1})=\|\hat{\pi}(x_{0:n-1})-\pi\|_{p}$.  Assume that $x=x_{0:n-1}$ and $\tilde{x}=\tilde{x}_{0:n-1}$ are two realizations of the random sequence $X_{0:n-1}$, which differ at a single element. By employing the reverse triangle inequality for the $\ell_p$-norm we obtain
\begin{align}\label{eq:Lipschitz_example_1}
\left|f(x)-f(\tilde{x})\right|&\leq \|\hat{\pi}(x_{0:n-1})-\hat{\pi}(\tilde{x}_{0:n-1})\|_{p}\nonumber\\&\leq \frac{\sqrt[p]{2}}{n}= \frac{\sqrt[p]{2}}{n} \tilde{d}_{1,\mathbb{E}^n}(x_{0:n-1},\tilde{x}_{0:n-1})=\frac{\sqrt[p]{2}}{n} \sum_{k=0}^{n-1}\mathbbm{1}(x_k\neq \tilde{x}_k).
\end{align} 
Clearly, (\ref{eq:Lipschitz_example_1}) implies that $\|f\|_{\rm Lip}=\sqrt[p]{2}/n$ for any $x_{0:n-1},\tilde{x}_{0:n-1}$ (not necessarily different at a single element). This can be easily seen by expressing $\hat{\pi}(x_{0:n-1})=(1/n)\sum_{k=0}^{n-1}\sum_{i=1}^{K}\mathbbm{1}(x_{k}=i)e_i$, where $\{e_1,\ldots, e_K\}$ is the canonical basis in $\mathbb{R}^K$. 

Consider the stationary chain $(\mathrm{P},\pi)$. An application of the theorem of Bobkov and G\"{o}tze (one-sided version) combined with Theorem \ref{thm:Marton} gives
\[
P_{\pi}(\|\hat{\pi}(X_{0:n-1})-\pi\|_{p}\geq E_{\pi}[\|\hat{\pi}(X_{0:n-1})-\pi\|_{p}]+\epsilon)\leq e^{-2^{1-2/p}n \epsilon^2 (1-r)^2}, \ \ \forall \epsilon>0
\]
or equivalently, $\forall \epsilon>E_{\pi}[\|\hat{\pi}(X_{0:n-1})-\pi\|_{p}]$,
\begin{equation*}
P_{\pi}(\|\hat{\pi}(X_{0:n-1})-\pi\|_{p}\geq \epsilon)\leq e^{-2^{1-2/p}n\left(\epsilon-E_{\pi}[\|\hat{\pi}(X_{0:n-1})-\pi\|_{p}]\right)^2 (1-r)^2}.
\end{equation*}
Theorem \ref{thm:Paulin} now implies that $\forall \epsilon>E_{\pi}[\|\hat{\pi}(X_{0:n-1})-\pi\|_{p}]$,
\begin{equation}
P_{\varrho}(\|\hat{\pi}(X_{0:n-1})-\pi\|_{p}\geq \epsilon)\leq \left\|\frac{\varrho}{\pi}\right\|_{2,\pi} e^{-2^{-2/p}n\left(\epsilon-E_{\pi}[\|\hat{\pi}(X_{0:n-1})-\pi\|_{p}]\right)^2 (1-r)^2}.
\end{equation}
Finally, for any $\delta\in (0,1)$ and any $\epsilon>E_{\pi}[\|\hat{\pi}(X_{0:n-1})-\pi\|_{p}]$, $P_{\varrho}(\|\hat{\pi}(X_{0:n-1})-\pi\|_{p}\geq \epsilon)\leq \delta$ for any $n$ such that 
\begin{equation}\label{eq:sample_complexity_Markov_1}
n\geq \frac{2^{\frac{2}{p}} \log\left(\frac{\left\|\frac{\varrho}{\pi}\right\|_{2,\pi}}{\delta}\right)}{\left(\epsilon-E_{\pi}[\|\hat{\pi}(X_{0:n-1})-\pi\|_{p}]\right)^2(1-r)^2}.
\end{equation}
$\newline$

\noindent\textbf{Approach 2: A union bound approach}

We may try to eliminate the problem of   $\epsilon$  being \emph{bounded away from zero} by using the observation that any $\ell_p$-norm is separable in the corresponding coordinates. We have
\begin{align*}
P_{\pi}(\|\hat{\pi}-\pi\|_{p}\geq \epsilon)&=P_{\pi}\left(\sum_{i=1}^K|\hat{\pi}_i-\pi_i|^p\geq \epsilon^p\right)\leq \sum_{i=1}^K P_{\pi}\left(|\hat{\pi}_i-E_{\pi}[\hat{\pi}_i]|\geq \frac{\epsilon}{\sqrt[p]{K}}\right),
\end{align*}
where the union bound and the fact that $\hat{\pi}_i$ are unbiased estimators $\forall i$ have been used. In this case $f(x)=f(x_{0:n-1})=\hat{\pi}_i(x_{0:n-1})$, therefore 
\[
|f(x)-f(\tilde{x})| \leq \frac{1}{n} \tilde{d}_{1,\mathbb{E}^n}(x_{0:n-1},\tilde{x}_{0:n-1}), \ \ \forall x_{0:n-1},\tilde{x}_{0:n-1}\in \mathbb{E}^n.
\]
By combining the theorem of Bobkov and G\"{o}tze with Theorems \ref{thm:Paulin} and \ref{thm:Marton} we obtain:
\begin{equation}\label{eq:scalarized_bound_Markov_1}
P_{\varrho}(\|\hat{\pi}(X_{0:n-1})-\pi\|_{p}\geq \epsilon)\leq \sqrt{2K}\left\|\frac{\varrho}{\pi}\right\|_{2,\pi} e^{-K^{-2/p}n\epsilon^2 (1-r)^2}.
\end{equation}

Finally, for any $\delta\in (0,1)$ and any $\epsilon>0$, $P_{\varrho}(\|\hat{\pi}(X_{0:n-1})-\pi\|_{p}\geq \epsilon)\leq \delta$ for any $n$ such that 
\begin{equation}\label{eq:sample_complexity_Markov_2}
n\geq \frac{K^{\frac{2}{p}} \log\left(\frac{\sqrt{2K}\left\|\frac{\varrho}{\pi}\right\|_{2,\pi}}{\delta}\right)}{\epsilon^2(1-r)^2}.
\end{equation}
\vspace{0.01cm}

\noindent\textbf{Approach 3: Application of Theorem \ref{thm:MainThm}}

 For simplicity, we will work with the first bound in (\ref{eq:MainThm_1}). Note that  by working with both bounds in (\ref{eq:MainThm_1}) we can only obtain an improvement of the derived sample complexity. 
 
 We observe that by the usual norm equivalence constants, the definition of $\tau_p$ in (\ref{eq:MainThm_3}) and the hierarchy of $\ell_p$-norms in $\mathbb{R}^k$ we have that $\tau_1\leq \sqrt{k}$, $\tau_2=1$, $\tau_p\leq 1$ for any $p>2$ and $\tau_{p_1}\leq \tau_{p_2}$ for $p_1\geq p_2$. By Theorems \ref{thm:MainThm}, \ref{thm:Paulin} and \ref{thm:Marton}  we obtain that $ \forall \epsilon>0$ and $\forall \varepsilon \in (0,1]$,
\begin{equation}\label{eq:Bobkov_Gotze_Extension_Markov_probability_bound}
P_{\varrho}(\|\hat{\pi}(X_{0:n-1})-\pi\|_{p}\geq \epsilon)\leq \left\|\frac{\varrho}{\pi}\right\|_{2,\pi} \left(1+\frac{2}{\varepsilon}\right)^{\frac{K}{2}}e^{-\frac{n\epsilon^2(1-\varepsilon)^2(1-r)^2}{2 \tau_p^2}}.
\end{equation}
Therefore, for any $\delta\in (0,1)$, any $\epsilon>0$ and any $\varepsilon\in (0,1]$, $P_{\varrho}(\|\hat{\pi}(X_{0:n-1})-\pi\|_{p}\geq \epsilon)\leq \delta$ for any $n$ such that 
\begin{equation}\label{eq:sample_complexity_Markov_3}
n\geq \frac{2\tau_p^2}{\epsilon^2(1-\varepsilon)^2(1-r)^2}\left[\frac{K}{2}\log\left(1+\frac{2}{\varepsilon}\right)+\log\left(\left\|\frac{\varrho}{\pi}\right\|_{2,\pi}\right)+\log\left(\frac{1}{\delta}\right)\right].
\end{equation}

\noindent\textbf{Sample Complexity Comparisons}

We first note that $\left\|\varrho/\pi\right\|_{2,\pi}\leq 1/\sqrt{\min_{i\in \mathbb{E}}\pi(i)}$ and often $\min_{i\in \mathbb{E}}\pi(i)\asymp 1/K^{m}$ for some $m\geq 1$. For a rough complexity comparison between (\ref{eq:sample_complexity_Markov_1}) and (\ref{eq:sample_complexity_Markov_3}) consider for simplicity the special case of i.i.d. random variables and $p=1$. In this setting, we correspondingly work with $P_{\pi}(\cdot)$ only (only $P_{\pi}(\cdot)$ is meaningful). It turns out that  $E_{\pi}[\|\hat{\pi}-\pi\|_1]\lesssim K/\sqrt{n}$. Then, (\ref{eq:sample_complexity_Markov_1}) and (\ref{eq:sample_complexity_Markov_3}) are orderwise the same, but without the problem of $\epsilon$ being bounded away from zero in (\ref{eq:sample_complexity_Markov_3}).  Further, (\ref{eq:sample_complexity_Markov_3}) is better by  a \textbf{logarithmic in $K$ factor} over (\ref{eq:sample_complexity_Markov_2}) for $p\in\{1,2\}$. For chains such that $\min_{i\in \mathbb{E}}\pi(i)\asymp 1/e^K$, (\ref{eq:sample_complexity_Markov_3}) is better by  a \textbf{$K$ factor} over (\ref{eq:sample_complexity_Markov_2}) for $p\in\{1,2\}$. We also note that depending on the geometry of $f(\mathbb{X})$, $f$ and $\nu$, $\tau_p$ may or may not have a favorable value for a particular $p$. More specifically, it is possible that the last approach is orderwise better than the union bound approach for some choices of $p$, primarily for $p\in \{1,2\}$, while it is worse for other values of $p$, depending on the particular problem at hand. 

\section*{Acknowledgment}

This work was supported by the ONR Grant Navy N00014-19-1-2566.


\bibliography{mybibfile_elsevier}

\begin{thebibliography}{18}
\providecommand{\natexlab}[1]{#1}
\providecommand{\url}[1]{\texttt{#1}}
\expandafter\ifx\csname urlstyle\endcsname\relax
  \providecommand{\doi}[1]{doi: #1}\else
  \providecommand{\doi}{doi: \begingroup \urlstyle{rm}\Url}\fi

\bibitem[Bobkov and G{\"o}tze(1999)]{bg99}
Sergej~G Bobkov and Friedrich G{\"o}tze.
\newblock Exponential integrability and transportation cost related to
  logarithmic sobolev inequalities.
\newblock \emph{Journal of Functional Analysis}, 163\penalty0 (1):\penalty0
  1--28, 1999.

\bibitem[Boucheron et~al.(2013)Boucheron, Lugosi, and Massart]{blm13}
St{\'e}phane Boucheron, G{\'a}bor Lugosi, and Pascal Massart.
\newblock \emph{Concentration inequalities: A nonasymptotic theory of
  independence}.
\newblock Oxford university press, 2013.

\bibitem[Br{\'e}maud(2013)]{b13}
Pierre Br{\'e}maud.
\newblock \emph{Markov chains: Gibbs fields, Monte Carlo simulation, and
  queues}, volume~31.
\newblock Springer Science \& Business Media, 2013.

\bibitem[Dembo and Zeitouni(1998)]{zd98}
Amir Dembo and Ofer Zeitouni.
\newblock \emph{Large Deviations Techniques and Applications}.
\newblock Springer-Verlag, New York, 1998.

\bibitem[Gromov(2007)]{g07}
Mikhail Gromov.
\newblock \emph{Metric structures for Riemannian and non-Riemannian spaces}.
\newblock Springer Science \& Business Media, 2007.

\bibitem[Hsu et~al.(2012)Hsu, Kakade, and Zhang]{hkz12}
Daniel Hsu, Sham Kakade, and Tong Zhang.
\newblock A tail inequality for quadratic forms of subgaussian random vectors.
\newblock \emph{Electronic Communications in Probability}, 17\penalty0
  (52):\penalty0 1--6, 2012.

\bibitem[Lattimore and Szepesv{\'a}ri(2020)]{lc20}
Tor Lattimore and Csaba Szepesv{\'a}ri.
\newblock \emph{Bandit algorithms}.
\newblock Cambridge University Press, 2020.

\bibitem[Levin and Peres(2017)]{lpw09}
David~A Levin and Yuval Peres.
\newblock \emph{Markov chains and mixing times}, volume 107.
\newblock American Mathematical Soc., 2017.

\bibitem[Marton(1986)]{m86}
Katalin Marton.
\newblock A simple proof of the blowing-up lemma (corresp.).
\newblock \emph{IEEE Transactions on Information Theory}, 32\penalty0
  (3):\penalty0 445--446, 1986.

\bibitem[Marton(1996)]{m96}
Katalin Marton.
\newblock Bounding $\bar{d}$-distance by informational divergence: A method to
  prove measure concentration.
\newblock \emph{The Annals of Probability}, 24\penalty0 (2):\penalty0 857--866,
  1996.

\bibitem[McDiarmid(1997)]{mc97}
Colin McDiarmid.
\newblock Centering sequences with bounded differences.
\newblock \emph{Combinatorics Probability and Computing}, 6\penalty0
  (1):\penalty0 79--86, 1997.

\bibitem[Paulin(2015)]{p15}
Daniel Paulin.
\newblock Concentration inequalities for {M}arkov chains by {M}arton couplings
  and spectral methods.
\newblock \emph{Electronic Journal of Probability}, 20\penalty0 (79):\penalty0
  1--32, 2015.

\bibitem[Raginsky and Sason(2014)]{rs13}
Maxim Raginsky and Igal Sason.
\newblock Concentration of measure inequalities in information theory,
  communications, and coding.
\newblock \emph{arXiv:1212.4663}, pages 1--180, 2014.

\bibitem[Rezakhanlou(2015)]{r15}
Fraydoun Rezakhanlou.
\newblock \emph{Lectures on the large deviation principle}.
\newblock Lecture Notes, Math UC Berkeley, 2015.

\bibitem[Van~Handel(2016)]{vh16}
Ramon Van~Handel.
\newblock \emph{Probability in High Dimension}.
\newblock APC 550 Lecture Notes, Princeton University, 2016.

\bibitem[Vershynin(2018)]{v18}
Roman Vershynin.
\newblock \emph{High-dimensional probability: An introduction with applications
  in data science}, volume~47.
\newblock Cambridge University Press, 2018.

\bibitem[Villani(2008)]{v08}
C{\'e}dric Villani.
\newblock \emph{Optimal transport: old and new}, volume 338.
\newblock Springer Science \& Business Media, 2008.

\bibitem[Wainwright(2019)]{w19}
Martin~J Wainwright.
\newblock \emph{High-dimensional statistics: A non-asymptotic viewpoint},
  volume~48.
\newblock Cambridge University Press, 2019.

\end{thebibliography}

\end{document}